\documentclass{amsart}
\usepackage{amsmath,amssymb}
\usepackage[UKenglish] {babel}

\usepackage{xy}
\xyoption{all}

\def\a{\alpha} \def\b{\beta} \def\d{\delta} \def\e{\epsilon} \def\f{\varphi}
\def\l{\lambda} \def\s{\sigma} \def\R{\mathbb{R}} \def\C{\mathbb{C}}
\def\H{\mathbb{H}} 
 \def\I{\mathbb{I}}
 \def\minus{\smallsetminus}
\def\({\left(} \def\){\right)} 
\def\<{\langle} \def\>{\rangle}
\def\w{\wedge}

\def\inv{^{-1}}

\renewcommand\ge{\geqslant}
\renewcommand\le{\leqslant}
\newcommand\ie{i.e.}

\newcommand\forget[1]{}

\newcommand\cs{\mathcal{S}}
\newcommand\cc{\mathcal{C}}

\DeclareMathOperator{\spann}{span}
\DeclareMathOperator{\rt}{R}
\DeclareMathOperator{\lt}{L}

\DeclareMathOperator{\Aut}{Aut}

\DeclareMathOperator{\so}{SO}

\DeclareMathOperator{\Mor}{Mor}
\DeclareMathOperator{\idp}{I}
\renewcommand{\Im}{\mathop{\rm Im}\nolimits}

\newenvironment{smatrix}{\left(\begin{smallmatrix}}{\end{smallmatrix}\right)}

\newcommand\pmatr[1]{\begin{pmatrix}#1\end{pmatrix}}

\newtheorem{thm}{Theorem}
\newtheorem{lma}[thm]{Lemma}
\newtheorem{prop}[thm]{Proposition}
\newtheorem{cor}[thm]{Corollary}
\theoremstyle{definition}

\theoremstyle{remark}
\newtheorem{rmk}[thm]{Remark}

\title[Four-dimensional power-commutative real division algebras]{Classification of the four-dimensional power-commutative real division algebras}
\author{Erik Darp\"{o} \and Abdellatif Rochdi}
\email{erik.darpo@maths.ox.ac.uk, abdellatifroc@hotmail.com}
\address{Erik Darp\"o: Mathematical Institute \\ 24-29 St Giles'\\
Oxford OX1~3LB\\ United Kingdom.
\linebreak
Abdellatif Rochdi: D\'epartement de Math\'ematiques et Informatique\\
Facult\'e des Sciences Ben~M'Sik\\ Universit\'e Hassan~II -- Mohammedia\\ B.P.~7955\\
Casablanca\\ Morocco 
}

\begin{document}
\selectlanguage{UKenglish}
\date{}

\begin{abstract}
  A classification of all four-dimensional power-commutative real division algebras is
  given. 

  It is shown that every four-dimensional power-commutative real division algebra is an
  isotope of a particular kind of a quadratic division algebra. The description of such
  isotopes in dimension four and eight is reduced to the description of quadratic division
  algebras. In dimension four this leads to a complete and irredundant
  classification. As a special case, the finite-dimensional power-commutative real division
  algebras that have a unique non-zero idempotent are characterised.
\end{abstract}

\maketitle

\section{Introduction}

An algebra is said to be \emph{power-commutative} if any subalgebra generated by a single
element is commutative. While this is always true for associative algebras, there are
plenty of non-associative algebras that fail to satisfy the power-commutative condition.
Indeed,
the power-commutative algebra structures on a
finite-dimensional vector space $V$ over a field $k$ form a subvariety of the variety
of all algebra structures on $V$, which is proper whenever $\dim V>1$.

On the other hand, the class of power-commutative algebras contains most types of
non-associative algebras that have been studied up to now. Commutative algebras and Lie
algebras are trivial examples. Some additional examples are:

\begin{enumerate}
\item
  \emph{Power-associative algebras:}  An algebra is power-associative if any element
  generates an associative subalgebra. This means that the powers $x^m$, $m\ge1$ are uniquely
  defined, and $x^mx^n=x^{m+n}$. Hence these algebras are power-commutative.
  All quadratic algebras and all alternative algebras are power-associative.
\item
  \emph{Flexible algebras:} This class is defined by the identity
  $(xy)x=x(yx)$. A subclass is the class of all non-commutative Jordan algebras, defined by
  the above equation and $x(yx^2)=(xy)x^2$.
  In a flexible algebra over a ground field of characteristic different from two, any
  commutative subset generates a commutative subalgebra \cite{raffin50}.
  Thus such an algebra is power-commutative.
\end{enumerate}

A \emph{division algebra} is an algebra $A$ in which the linear maps $\lt_a:A\to
A,\:x\mapsto ax$ and $\rt_a:A\to A,\:x\mapsto xa$ are bijective for all non-zero 
$a\in A$. An equivalent condition for finite-dimensional algebras is that $xy=0$ only if
$x=0$ or $y=0$.
Over the real numbers $\R$, every finite-dimensional associative division algebra is
isomorphic to $\R$, $\C$ or the quaternions $\H$ (Frobenius \cite{frobenius78} 1878).
Replacing associativity by the weaker condition of alternativity\footnote{An
  algebra is called \emph{alternative} if every subalgebra generated by two elements is
  associative.} gives only one additional isomorphism class, the one of the octonion
algebra $\mathbb{O}$ (Zorn \cite{zorn31} 1931). 
Other classical theorems assert that the dimension of a real division algebra, if finite,
is one, two, four or eight (Bott and Milnor \cite{bottmilnor}, Kervaire
\cite{kervaire58} 1958), and that every commutative 
finite-dimensional real division algebra has dimension one or two (Hopf \cite{hopf} 1940).
The last result plays a decisive role in the present article, as it limits the dimension
of a subalgebra generated by a single element to a most two.

Much of the theory on real division algebras that has been developed after 1958 emanates
from a series of works by Benkart, Britten and Osborn:
\cite{osborn62,isotopy,bo81a,bo81b,bbo82}.

In the pioneering article \cite{osborn62}, Osborn initiated the study of quadratic
division algebras. An algebra $A$ is called \emph{quadratic} if it has an identity
element, and $1,x,x^2$ are linearly dependent for all $x\in A$. Dieterich \cite{zur}
finished the classification of all quadratic division algebras of dimension four over
$\R$. The problem of classifying all eight-dimensional quadratic real division algebras is
still open.

Flexible finite-dimensional real division algebras were studied by Benkart, Britten and
Osborn in \cite{bbo82}. Based on their approach, a classification has been accomplished in
recent years, through the contributions \cite{malaga,coll,nform}.

Several independent classifications of the two-dimensional real division algebras exist:
\cite{burdujan,gottschling,hupe04,dieterich05}.
Whereas the former two employ the classical approach with multiplication tables,
\cite{hupe04} and \cite{dieterich05} build on the notion of isotopy, which was introduced
by Albert \cite{albert42a} and used by Benkart, Britten and Osborn in \cite{isotopy}.
The concept of isotopy plays a key role also in the present article, as will explained
later in this section.

Some other important contributions to the theory of real division algebras are the
treatment of division algebras via their derivation algebras
(\cite{bo81a,bo81b,rochdi95,dz04b,perez06}) and the theory of finite-dimensional absolute
valued algebras (its development until 2004 is surveyed in \cite{rodriguez04}, some more
recent contributions are \cite{roro09,ava,ckmmrr09}). 

In spite of these advances, the problem of classifying all finite-dimensional
power-commutative real division algebras has remained untackled.
The present paper is devoted to the solution of this problem in dimension four.
A pivotal result is Theorem~\ref{4dim}, formulated below. It asserts that every
four-dimensional power-commutative real division algebra $A$ possesses a non-zero
idempotent $e$ that is contained in every two-dimensional subalgebra of $A$. Such an
idempotent we shall call \emph{omnipresent}.
Section~\ref{sec:lemmata} contains some general, technical results, valid for any
finite-dimensional power-commutative division algebra over $\R$. 
Section~\ref{sec:eq} investigates the structure of power-commutative real division
algebras containing an omnipresent idempotent, reaching a classification in the
four-dimensional case.
As a special case, Section~\ref{sec:enidp} treats real division algebras having a unique
non-zero idempotent. 
Finally, the proof of Theorem~\ref{4dim} is given in Section~\ref{sec:4dim}.

All algebras and vector spaces occurring in this article are finite-dimensional over the
real numbers. In the sequel, these assumptions are implicit when we use the words
'algebra', 'vector space' etc.

Hopf's theorem on commutative division algebras allow for a useful, alternative
characterisation of power-commutativity: 
a division algebra $A$ is power-com\-mu\-ta\-tive if and only if for all $x\in A$ the
subspace $\spann\{x,x^2\}\subset A$ is a subalgebra, and the identity
\begin{equation}
  \label{3pass}
  xx^2=x^2x
\end{equation}
holds.
An algebra satisfying Equation~(\ref{3pass}) is called \emph{third-power associative}.

\begin{thm} \label{4dim}
  Every four-dimensional power-commutative real division algebra contains an omnipresent
  idempotent.
\end{thm}
Clearly, a power-commutative real division algebra $A$ of dimension greater than two can
contain at most one omnipresent idempotent.
Note also that a non-zero idempotent $e\in A$ is omnipresent if and only if
$\spann\{e,x\}$ is a subalgebra of $A$ for all $x\in A$: if $e$ is omnipresent and $e,\,x$
are non-proportional then, by Lemma~\ref{sdzd}(2) in Section~\ref{sec:lemmata},
$\spann\{e,x\}$ contains a non-idempotent $y$, the subalgebra generated by which
coincides with $\spann\{e,y\}=\spann\{e,x\}$. 
In particular, an omnipresent idempotent commutes with all elements in $A$.

Given an algebra $A$ and two invertible linear maps $S,T:A\to A$, a new algebra structure
on $A$ is defined by $x\circ y= (Sx)(Ty)$.
The algebra $A_{ST}=(A,\circ)$ obtained in this way is called the
\emph{(principal) isotope of $A$ given by $S$ and $T$}.
Since $S$ and $T$ are invertible, the algebra $A_{ST}$ is a division algebra if and
only if $A$ is a division algebra.
In this article we will be interested in isotopes $A_{ST}$ where $S=T$, and we write
$A_T=A_{TT}$. If $A$ is a division algebra and $e\in A$ a non-zero idempotent commuting
with all elements in $A$, then $A_{\lt_e^{-1}}$ is unital with identity element $e$.

In any quadratic algebra $B$, the set
$$\Im B = \{x\in B\minus \R1_B \mid x^2\in \R1_B\}\cup\{0\} \:\subset\: B \,.$$
is a subspace of $B$, and $B=\R1_B\oplus \Im B$
(Frobenius \cite{frobenius78} 1878). 
Let $B$ be a quadratic division algebra. A linear map $T:B\to B$ shall be said to be
\emph{planar} if it is invertible, $T(1)=1$, and every subalgebra of $B$ is
invariant under $T$.
For any $x\in B$, the subspace $\spann\{1,x\}\subset B$ is a subalgebra, so
$T(x)\in\spann\{1,x\}$. On the other hand, since the identity element is contained
in every subalgebra of $B$, an invertible linear map $T:B\to B$ satisfying $T(1)=1$ and
$T(x)\in\spann\{1,x\}$ for all $x\in B$ is planar.

Let $P:B\to \Im B$ be the projection onto $\Im B$ along $\R1$. 
If $T:B\to B$ is a planar map, the linear map $T'=PT:\Im B\to \Im B$ satisfies
$T'(\R v)=\R v$ for all $v\in\Im B$. Hence $T'=\l\I_{\Im B}$ for some scalar
$\l\in\R\minus\{0\}$, and $T(v)=\s(v)1+\l v$ for all $v\in\Im B$, where $\s:\Im B\to\R$ is
a linear form, determined by $T$.
Conversely, $T$ is given entirely by the form $\s$ and the scalar $\l$.
In a basis $\underline{v}=(1,v_2,\ldots,v_n)$ of $B$, with $v_i\in\Im B$, the
matrix of $T$ takes the form
$$
  [T]_{\underline{v}}= \left(
  \begin{array}{c|ccc}
    1 & s_2&\cdots&s_n \\ \hline 
    0 & \l \\
    \vdots && \ddots \\
    0 &&& \l
  \end{array}
  \right).
$$

We shall say that $B_T$ is a \emph{planar isotope} of $B$ if $T:B\to B$ is planar. In
Section~\ref{sec:eq} we show that every power-commutative division algebra containing an
omnipresent idempotent is a planar isotope of a quadratic division algebra. 
In particular, in view of Theorem~\ref{4dim}, every four-dimensional power-commutative
division algebra is of this form. 

The following notation and conventions are used throughout the article.
In any algebra $A$, we define $[x,y]=xy-yx$ and $x\bullet y =xy+yx$.
The \emph{centraliser} of $x\in A$ is defined as 
$C(x)=\{y\in A \mid [x,y]=0\}$. 
An element $a\in A$ is said to be \emph{central} if $C(a)=A$.
The subalgebra of $A$ generated by a subset $M\subset A$ is denoted by $A(M)$, and
$\idp(A)$ denotes the set of non-zero idempotents in $A$.
The vector space $\R^m$ is understood to consist of column vectors, and linear maps
$\R^m\to\R^n$ are identified with $n\!\times\! m$-matrices in the natural way.
If $\mathcal{A}$ is a category in which a concept of dimension is defined,
$\mathcal{A}_{m}$ and $\mathcal{A}_{\ge m}$ denote the full subcategories of $\mathcal{A}$
formed by objects of dimension $m$ and objects of dimension greater than or equal to $m$,
respectively.

\section{Lemmata} \label{sec:lemmata}

In this section, we give some basic results on power-commutative real division
algebras.

Substituting $x+\l y$ for $x$ in the third-power associative identity (\ref{3pass}), we
get the equation 
$$\l\left( [x\bullet y,x] + [x^2,y] \right) +
 \l^2\left( [y^2,x] + [x\bullet y,y] \right)=0$$
valid for arbitrary $x$ and $y$, and all $\l\in\R$. This is possible only if each of the
terms in the sum equals zero, that is if
\begin{equation}
  \label{kvadkomm}
  [x\bullet y , y] = [x,y^2]
\end{equation}
for all $x$ and $y$.
In particular, if $e=y$ is an idempotent the equation
\begin{equation}
  \label{idpkomm}
  [e,e\bullet x - x]=0
\end{equation}
holds. If $e$ and $x$ commute, then (\ref{idpkomm}) reduces to $[e,2ex]=0$. Thus we have
the following result.
\begin{lma}\label{Leinvariant}
  In a third-power associative algebra, the centraliser $C(e)$ of any idempotent $e$ is
  invariant under $\lt_e$ (and hence also under $\rt_e$).
\end{lma}

For any two idempotents $e,f$ in an arbitrary algebra $A$,
$(e+f)^2=e+f +e\bullet f$ and $(e-f)^2=e+f -e\bullet f=2(e+f) -(e+f)^2$ hold. 
Thus $(e-f)^2\in A(e+f)\cap A(e-f)$.
Suppose $A$ is a power-commutative division algebra and $[e,f]\ne0$. 
Then $e-f\not\in A(e+f)$, because otherwise $A(e,f)=A(e+f,e-f)=A(e+f)$, and 
$A(e,f)$ would be commutative, implying $[e,f]=0$.
As any subalgebra of $A$ generated by a single element has dimension at most two,
$A(e+f)\cap A(e-f)=\spann\{(e-f)^2\}$ and consequently $\left((e-f)^2\right)^2=\mu(e-f)^2$
for some $\mu\in\R\minus\{0\}$. Moreover $[(e-f)^2,e-f]=0=[(e-f)^2,e+f]$, but
$[e,f]\ne0$, so $(e-f)^2\not\in\spann\{e+f \,,\,e-f\}=\spann\{e,f\}$.
 This proves the following lemma: 
\begin{lma} \label{eminusfkvadrat}
  If $e,f$ are non-commuting idempotents in a power-commutative division algebra $A$, then
  $(e-f)^2$ is a scalar multiple of an idempotent, not contained in $\spann\{e,f\}$.
\end{lma}

Our next lemma concerns the case of commuting idempotents.

\begin{lma} \label{kommidpsubalg}
  Let $A$ be a power-commutative division algebra. Then every set of mutually commuting
  idempotents spans a (commutative) subalgebra of $A$. Thus no more than two linearly
  independent idempotents can commute with each other.
\end{lma}

\begin{proof}
  It suffices to show that if $e,f\in A$ are two distinct idempotents commuting with each
  other, then $\spann\{e,f\}\subset A$ is a subalgebra.
  The subalgebra $A(e-f)\subset A$ is a division algebra, and so contains an element $x\in
  A(e-f)$ such that $x(e-f)=e-f$. But $(e+f)(e-f)=e-f$, hence $((e+f)-x)(e-f)=0$, which
  implies $e+f=x\in A(e-f)$. Thus $\spann\{e,f\}=\spann\{e+f\,,\,e-f\}= A(e-f)$. 
\end{proof}

Finally, we recall three known results in the theory of division algebras, which will
be used in the sequel.

\begin{lma} \label{sdzd}
  \begin{enumerate}
  \item \cite[Theorem~1]{segre} Every algebra in which $x^2=0$ only if $x=0$
    contains a non-zero idempotent.
  \item \cite[Proposition~2.2]{dz04b} Any two-dimensional subspace of an algebra satisfying
    the condition in (1) contains at most three non-zero idempotents.
  \item \cite[Lemma~2.3]{coll} Given an idempotent $e$ in a two-dimensional
    commutative division algebra $A$, there exists an element $x\in A$ such that
    $x^2=-e$.
  \end{enumerate}
\end{lma}

\section{Algebras with an omnipresent idempotent} \label{sec:eq}

In this section, we prove that all power-commutative division algebras with an omnipresent
idempotent are planar isotopes of quadratic division algebras
(Theorem~\ref{eqisotope}), and describe the morphisms between such algebras in
dimension four and eight (Proposition~\ref{morphisms}).
From this a complete and irredundant classification in dimension four is deduced.
In dimension eight, further progress depends on the understanding of the quadratic case, of
which still rather little is known. 

We begin by listing some criteria
for a power-commutative division algebra to be quadratic and to have an omnipresent
idempotent, respectively. The items (\ref{q}) and (\ref{eq}) are not necessary for any of
the other results in this article, they are included here only for their own interest.
Item~(\ref{Icomm}), in contrast, plays an important role in the proof of Theorem~\ref{4dim}.

\begin{prop} \label{eqvillkor}
  Let $A$ be a power-commutative division algebra.
  \begin{enumerate}
  \item The algebra $A$ is quadratic if and only if it contains an identity
    element. \label{q} 
  \item The algebra $A$ contains an omnipresent idempotent if and only if any
    two-dimensional subalgebras $B_1,B_2\subset A$ have a non-trivial
    intersection. \label{eq}
  \item If $\idp(A)$ is commutative, then $A$ contains an omnipresent
    idempotent. \label{Icomm}
  \end{enumerate}
\end{prop}

\begin{proof}
\ref{q}) A quadratic division algebra is a unital division algebra $A$ in which 
$\dim A(x)\le2$ for all $x\in A$.
Thus a power-commutative division algebra, by default satisfying the latter condition, is
quadratic if and only if it possesses and identity element.

\ref{eq})
The implication $"\Rightarrow"$ is immediate from the definition of an omnipresent
idempotent.
For the other direction, we must prove that all
two-dimensional subalgebras have a common intersection, which then
contains the omnipresent idempotent $e$. For this, it suffices to show
that if $B_1, B_2, B_3\subset A$ are such subalgebras, then
$B_1\cap B_2\cap B_3\neq 0.$ 
By assumption, $B_i\cap B_j\neq 0,$ and there is a unique non-zero
idempotent $e_{ij}$ in $B_i\cap B_j.$ Since $e_{12}, e_{13}\in B_1,$ we have 
$[e_{12},e_{13}]=0$ and analogously, $[e_{12}, e_{23}]=[e_{13}, e_{23}]=0.$ Thus, by
Lemma~\ref{kommidpsubalg}, the idempotents $e_{ij}$ span a commutative subalgebra of $A$,
and hence they are linearly dependent. If $e_{12}, e_{13}, e_{23}$ are
distinct, then $e_{23}\in \spann\{e_{12}, e_{13}\}$, and $B_1=B_2=B_3,$. 
If two of them are equal, say $e_{12}=e_{13}$, then $e_{12}=e_{13}\in B_1\cap
B_2\cap B_3$.

\ref{Icomm})
The result is trivial in dimension one and two. Thus assume $\dim A\ge4$.
By Lemma~\ref{kommidpsubalg}, $B=\spann\idp(A)$ is a commutative subalgebra of $A$, hence
$\dim B\le2$ and Lemma~\ref{sdzd}(2) gives $|\idp(A)|\le3$.
Now for any $x\in A\minus B$, the subalgebra $A(x)\subset A$ intersects
$B$ in a one-dimensional subspace, containing a unique non-zero idempotent
$e_x$.
This defines a map
$$\iota:A\!\minus B\,\to\, \idp(A),\;x\mapsto e_x.$$

Given a vector space $V$ and $m\in\mathbb{N}$, let $g(V,m)$ denote the Grassmanian
manifold of $m$-dimensional subspaces of $V$.
Now the map $\iota$ factors as $\iota=\iota_3\iota_2\iota_1$, where
\begin{align*}
  \iota_1&:A\minus B\to g(A,2),\:x\mapsto A(x),\\ 
  \iota_2&:\iota_1(A\minus B)\to g(B,1),\:U\mapsto B\cap U,  \\
  \iota_3&:\iota_2\iota_1(A\minus B)\to\idp(A),\:\spann\{e\}\mapsto e \,.
\end{align*}
Each of the maps $\iota_i$ are clearly continuous, and hence $\iota$ is continuous.
Since $B\subset A$ has codimension at least two, $A\!\minus B$ is
connected. Then $\iota(A\!\minus B)\subset\idp(A)$ is also connected, and thus a
singleton: $\iota(A\!\minus B) = \{e\}$. This means that $x^2\in\spann\{e,x\}$
for all $x\in A$, so $e$ is omnipresent.
\end{proof}

In particular, the third part of Proposition~\ref{eqvillkor} implies that if $e$ is the
only non-zero idempotent in a power-commutative division algebra $A$, then $e$ is
omnipresent in $A$ (indeed, this conclusion follows from Lemma~\ref{sdzd}(1) as well).

We now turn to the problem of describing all power-commutative division algebras with an
omnipresent idempotent. In Theorem~\ref{eqisotope} we will show that they are precisely
the planar isotopes of quadratic division algebras. The foundation is set by the following
lemma, which shows that every  power-commutative division algebra with an omnipresent
idempotent is isotopic to a quadratic division algebra.

\begin{lma} \label{eqlemma}
  Let $A$ be a power-commutative division algebra, and $e$ an omnipresent idempotent in
  $A$.
  \begin{enumerate}
  \item A subspace $U\subset A$ with $\dim U\ge2$ is a subalgebra of $A$ if and only if it
    is a subalgebra of  $A_{\lt_e\inv}$.
  \item The algebra $A_{\lt_e\inv}$ is quadratic, with identity element $e$. 
  \end{enumerate}
\end{lma}

\begin{proof}
  Every subalgebra $B$ of $A$ of dimension greater than one contains $e$ and,
  consequently, is invariant under $\lt_e$ and $\lt_e\inv$. In particular,
  $\lt_e\inv(x)\in\spann\{e,x\}$ for all $x\in A$. 
  Denoting by $\circ$ the multiplication in $A_{\lt_e\inv}$, we have 
  $x\circ y = \lt_e\inv(x)\lt_e\inv(y)\in B$ for all 
  $x,y\in B$, hence $B$ is a subalgebra of $A_{\lt_e\inv}$.
  
  Since $e$ is a central idempotent in $A$, it is an identity element in
  $A_{\lt_e\inv}$. Moreover, 
  $x\circ x=\left(\lt_e\inv(x)\right)^2\in\spann\{e,\lt_e\inv(x)\}=\spann\{e,x\}$, which
  proves that $A_{\lt_e\inv}$ is quadratic.

  Every subalgebra of a quadratic algebra contains the identity element. Hence, any
  subalgebra $B$ of $A_{\lt_e\inv}$ is invariant under $\lt_e$. By the same argument as
  above, it follows that $B$ is a subalgebra of $\left(A_{\lt_e\inv}\right)_{\lt_e}=A$.
\end{proof}

\begin{thm} \label{eqisotope}
\begin{enumerate}
\item An isotope $B_T$ of a quadratic division algebra $B$ of dimension at least four is
  power-commutative if and only if $T=\l S$ for some planar map $S:B\to B$ and non-zero
  number $\l\in\R$. 
  In this case, $1_B$ is an omnipresent idempotent in $B_T$.
\item Every power-commutative division algebra with an omnipresent idempotent is a planar
  isotope of a quadratic division algebra. 
\end{enumerate}

\end{thm}

\begin{proof}
1.
It is readily verified that $A=B_T=(B,\circ)$ is power-commutative, with omnipresent
idempotent $1_B$ if $T$ is a scalar multiple of a planar map. 

For the converse, let $[\,]_A$ be the commutator in $A$, whereas $[\,]$ denotes the
commutator in $B$.
Power-commutativity of $A$ gives  
$$0=[x\circ x,x]_A = \left[(Tx)^2,x \right] = \left[T\left((Tx)^2\right), Tx\right]$$
for all $x\in B$, equivalently, $\left[Tx^2, x\right]=0$.

Specifying $x\in\Im B$ in the last identity yields $[T(1_B),x]=0$ for all $x\in\Im B$, that
is, $T(1_B)=\l1_B\in\R 1_B$.
On the other hand, if $x=1_B+v$ with $v\in\Im B$, then $x^2=1_B+v^2+2v=\mu1_B+2v$,
$\mu\in\R$, and $[Tx^2,x]=[\mu T(1_B)+ 2Tv, 1_B+v]=2[Tv,v]$.
In any quadratic algebra, the centraliser of an element $y\not\in\R1$ is 
$C(y)=\spann\{1,y\}$, so $2[Tv,v]=[Tx^2,x]=0$ implies $Tv\in\spann\{1_B,v\}$.

Hence $Tx\in\spann\{1_B,x\}$ for any $x\in B$. It follows that $\frac{1}{\l}T$ is planar. 

2.
This is trivial in dimension one and two. In dimension greater than two, it follows from
the first part of the theorem together with Lemma~\ref{eqlemma}(2).
\end{proof}

The remainder of this section is devoted to the problem of describing isomorphism classes
of planar isotopes. In dimension four, this results in a complete and irredundant
classification. In dimension eight we obtain a general description, which essentially
reduces the classification problem of planar isotopes to the problem of classifying
dissident maps in seven-dimensional Euclidean space. The latter problem is still far from
being solved, though some progress has been made in recent years; see
\cite{dili,lindberg,liftings,diru10}.

\begin{prop} \label{morphisms}
  Let $A$ and $B$ be quadratic division algebras of dimension at least four, and $S$
  and $T$ planar maps on $A$ and $B$ respectively. Then
  $\Mor(A_S,B_T)=\{\f\in\Mor(A,B) \mid \f S = T\f \}$.
\end{prop}

\begin{proof}
We use $\circ$ to denote the multiplication in the isotopes $A_S$ and $B_T$.

It is easy to verify that any morphism $\f:A\to B$ satisfying $\f S = T\f$ is a
morphism of the isotopes $A_S\to B_T$.
For the converse, let $\f\in\Mor(A_S,B_T)$. Now $\f(1_A)\in B$ is an idempotent commuting
with $\f(A_S)\subset B$, which is a subalgebra of dimension at least four. Since the
centraliser of any element in $B$ outside $\R1_B$ has dimension two, it follows that
$\f(1_A)=1_B$.

For arbitrary $x\in A$, 
$$ \f S(x)=\f(1_A\circ x)= \f(1_A)\circ\f(x)= 1_B\circ\f(x)= T\f(x), $$
that is, $\f S = T\f$.
Moreover,
\begin{align*}
  \f(x\circ y) &= \f((Sx)(Sy)) \quad\mbox{and} \\
  \f(x)\circ\f(y) &= (T\f(x))(T\f(y)) = \f(Sx)\,\f(Sy)
\end{align*}
for all $x,y\in A$. Since $S:A\to A$ is bijective, this implies that $\f(xy)=\f(x)\f(y)$
for all $x,y\in A$, \ie, that $\f$ is a morphism of the quadratic algebras $A$ and
$B$.
\end{proof}

Suppose $A$ and $B$ are quadratic division algebras. Let $S:A\to A$ be the planar map
determined by linear form $\s:\Im A\to\R$ and scalar $\l\in\R\minus\{0\}$, and $T:B\to B$
the planar map determined by $\tau:\Im B\to\R$ and $\mu\in\R\minus\{0\}$. 
Then a morphism $\f:A\to B$
satisfies $\f S = T\f$ if and only if $\l=\mu$ and $\s=\tau\f$.  

The concept of a \emph{dissident triple} was introduced in \cite{dissalg}, as a way to
describe quadratic division algebras. As we shall see, the concept can be modified to
give a complete description not only of quadratic division algebras, but also of their
planar isotopes. First, we recapitulate some basic facts about dissident triples.

A \emph{dissident map} on a Euclidean space $V=(V,\<\,\>)$ is a linear map $\eta:V\w V\to V$
with the property that $u,v,\eta(u\w v)$ are linearly independent whenever $u,v\in V$ are.
A \emph{dissident triple} is a triple $(V,\xi,\eta)$ consisting of a Euclidean space $V$, a
linear map $\xi:V\w V\to\R$ and a dissident map $\eta:V\w V \to V$. Declaring as morphisms
$(V,\xi,\eta)\to(V',\xi',\eta')$ those orthogonal maps $\f:V\to V'$ that satisfy
$\xi=\xi'(\f\w\f)$ and $\f\eta=\eta'(\f\w\f)$, the class of dissident triples is endowed
with the structure of a category, which we denote by $\tilde{\mathcal{C}}$.

The relevance of dissident triples lies in the fact that they correspond precisely to
quadratic division algebras. Given a dissident triple $(V,\xi,\eta)$, the
algebra $\tilde{\mathcal{F}}(V,\xi,\eta)=\R\times V$ with multiplication defined by 
$$\pmatr{\a\\u}\pmatr{\b\\v} = \pmatr{\a\b -\<u,v\>+\xi(u\w v) \\ \a v+\b u +\eta(u\w v)}$$
is a quadratic division algebra. Moreover, defining $\tilde{\mathcal{F}}(\f)=\I_\R\times\f$ for
morphisms, $\tilde{\mathcal{F}}$ becomes an equivalence between the respective categories of dissident
triples and quadratic division algebras \cite{osborn62,revisited}.
For convenience, we shall identify the imaginary space
$\Im\tilde{\mathcal{F}}(V,\xi,\eta)=\{0\}\times V$ with $V$.

Set $\mathcal{T}=\{d\in\R^3 \mid 0<d_1\le d_2\le d_3\}$. Any triple
$(x,y,d)\in\R^3\times\R^3\times\mathcal{T}$ gives rise to a dissident triple
$\tilde{\mathcal{G}}(x,y,d)=(\R^3,\xi_x,\eta_{yd})$ on setting $\xi_x(u\w v)=\<\e_x(u),v\>$ and
$\eta_{yd}(u\w v)=(\e_y+\d_d)\pi(u\w v)$, where
\begin{equation*}
  \d_{d}=\begin{smatrix}
    d_1&0&0 \\
    0&d_2&0 \\
    0&0&d_3
  \end{smatrix},
\qquad 
  \e_x=\begin{smatrix}
    0&-x_3&x_2 \\
    x_3&0&-x_1 \\
    -x_2&x_1&0
  \end{smatrix}
\end{equation*}
and $\pi:\R^3\w\R^3\to\R^3$ is the standard vector product.
The set $\R^3\times\R^3\times\mathcal{T}$ is given the structure of a category by defining
a morphism $(x,y,d)\to (x',y',d')$ to be a map $\f\in\so(\R^3)$ satisfying
$(\f(x),\f(y),\f\d_d\f\inv)=(x',y',\d_{d'})$. We shall denote this category by
$\tilde{\mathcal{K}}$. The map $\tilde{\mathcal{G}}$ described above defines an
equivalence $\tilde{\mathcal{K}} \to \tilde{\mathcal{C}}_3$, acting on morphisms
identically \cite{ohman}.

By adding to a dissident triple $(V,\xi,\eta)$ a linear form $\s:V\to\R$ and a non-zero
scalar $\l$, the information corresponding to a planar map is encoded. 
Thus, we consider a category $\mathcal{C}$ whose objects are tuples $(V,\xi,\eta,\s,\l)$
where $(V,\xi,\eta)$ is a dissident triple, $\s$ a linear form on $V$, and
$\l\in\R\minus\{0\}$.
For $X=(V,\xi,\eta,\s,\l)$ and $X'=(V',\xi',\eta',\s',\l')$ in $\mathcal{C}$, we define
morphisms as
$$\Mor_{\mathcal{C}}(X,X') = 
\begin{cases}
   \{\f\in\Mor_{\tilde{\mathcal{C}}}((V,\xi,\eta),(V',\xi',\eta')) \mid \s=\f\s' \} &\text{if } \l=\l', \\
   \emptyset & \text{otherwise}.
 \end{cases}
 $$

Objects in $\mathcal{C}$ can be seen as diagrams in  
the category of vector spaces of the form \\
\centerline{
\xymatrix{
V\w V \ar[dr]_\xi \ar[rr]^\eta && V \ar[dl]^\s \\
& \R \ar@(dl,dr)_\l
} 
}
where $V$ is a Euclidean space, $\eta$ is dissident and $\l$ non-zero. Morphisms in
$\mathcal{C}$ then are diagram morphisms of the form $(\f\w\f,\f,\I_\R)$, with $\f$
orthogonal.

Let $\mathcal{D}$ be the category of power-commutative division algebras containing an
omnipresent idempotent. 
A functor $\mathcal{F}:\mathcal{C}\to\mathcal{D}$ is defined by
$\mathcal{F}(V,\xi,\eta,\s,\l)=\tilde{\mathcal{F}}(V,\xi,\eta)_{S}$, where $S$ is the
planar map determined by $\s$ and $\l$, and $\mathcal{F}(\f)=(\I_\R\times\f)$.

\begin{prop} \label{dissekv}
  The functor $\mathcal{F}:\mathcal{C}\to\mathcal{D}$ is faithful and dense. It induces an equivalence
  $\mathcal{C}_{\ge3}\to\mathcal{D}_{\ge4}$.
\end{prop}

\begin{proof}
  Faithfulness is immediate from the construction. Since every quadratic division algebra
  is isomorphic to $\tilde{\mathcal{F}}(V,\xi,\eta)$ for some dissident triple $(V,\xi,\eta)$,
  every planar isotope $A_T$ of a quadratic division algebra $A$ is isomorphic to some
  $\mathcal{F}(V,\xi,\eta,\s,\l)$. By Proposition~\ref{eqisotope}, these are precisely the
  power-com\-mu\-ta\-tive division algebras containing an omnipresent idempotent.

  Let $A_S=\mathcal{F}(V,\xi,\eta,\s,\l)$ and $B_T=\mathcal{F}(V',\xi',\eta',\s',\l')$,
  with $\dim A_S,\dim B_T\ge4$.   
  If $\psi:A_S\to B_T$ is a morphism then, by Proposition~\ref{morphisms},
  $\psi\in\Mor(A,B)$, $\l=\l'$ and $\s=\s'\psi|_{\Im A}$. Now $\psi\in\Mor(A,B)$ implies
  that $\psi=\tilde{\mathcal{F}}(\f)=\I_\R\times\f$ for some
  $\f:(V,\xi,\eta)\to(V',\xi',\eta')$. Furthermore, since
  $\l=\l'$ and $\s=\s'\psi|_{\Im A}=\s'\psi|_V=\s'\f$, it follows that $\f$ is a morphism in
  $\mathcal{C}$. Hence $\psi=\mathcal{F}(\f)$, and the induced functor
  $\mathcal{F}:\mathcal{C}_{\ge3}\to\mathcal{D}_{\ge4}$ is full, thus an equivalence.
\end{proof}

\begin{rmk}
  The multiplication in the algebra $\mathcal{F}(V,\xi,\eta,\s,\l)=(\R\times V\,,\,\circ)$
  is given by
  \begin{equation} \label{kvadmult}
    \pmatr{\a\\u}\circ\pmatr{\b\\v} = 
  \begin{pmatrix}
    \a\b +\a\s(v)+\b\s(u)+\s(u)\s(v) + \l^2\xi(u\w v) - \l^2\<u,v\> \\
    \l(\a+\s(u))v + \l(\b+\s(v))u + \l^2\eta(u\w v)
  \end{pmatrix} .
\end{equation}
\end{rmk}

Just as the category $\tilde{\mathcal{C}}$ was modified to correspond to all $e$-quadratic
power-commutative division algebras, $\tilde{\mathcal{K}}$ can be extended
to match the category $\mathcal{C}_3$. 
A Euclidean space $V$ identifies with its dual space $V^\ast$ via the map $v\mapsto\s_v$, $\s_v(u)=\<u,v\>$.
If $\f:(V,\xi,\eta)\to(V',\xi',\eta')$ is an isomorphism in $\tilde{\mathcal{C}}$ and
$u,v\in V$, then $\s_u=\s_v\f$ if and only if $u=\f^\ast(v)$, that is, since $\f$ is
orthogonal, if and only if $v=\f(u)$. In this case, $\f$ defines an isomorphism
$(V,\xi,\eta,\s_u,\l)\to(V',\xi',\eta',\s_v,\l)$ in $\mathcal{C}$.

Let $\mathcal{K}$ be the category with object class
$\R^3\times\R^3\times\R^3\times\mathcal{T}\times(\R\minus\{0\})$ and 
$$\Mor(\kappa,\kappa')= \{ \f\in\so(\R^3) \mid 
(\f(x),\f(y),\f(z),\f \d_d\f\inv,\l)=(x',y',z',\d_{d'},\l') \}$$
for $\kappa=(x,y,z,d,\l)$, $\kappa'=(x',y',z',d',\l')$ in $\mathcal{K}$.
Now, in view of the equivalence
$\tilde{\mathcal{G}}:\tilde{\mathcal{K}}\to\tilde{\mathcal{C}}_3$ and the above
consideration, the following result is readily verified.

\begin{prop} \label{konfekv}
  The functor $\mathcal{G}:\mathcal{K}\to\mathcal{C}_3$ defined by
  $\mathcal{G}(x,y,z,d,\l)=(\R^3,\xi_x,\eta_{yd},\linebreak \s_z,\l)$ and
  $\mathcal{G}(\f)=\f$ is an equivalence of categories.
\end{prop}

In order to classify all four-dimensional  power-commutative division algebras having an
omnipresent idempotent, it remains to determine a cross-section for the isomorphism
classes in $\mathcal{K}$, \ie, a set of objects $\cs\subset \mathcal{K}$ containing a
unique representative of each isomorphism class. By Propositions~\ref{dissekv} and
\ref{konfekv}, 
$\mathcal{F}\mathcal{G}(\cs)$ then is a cross-section for the isomorphism classes
in $\mathcal{D}_4$. The problem is similar to the classification problem for
$\tilde{\mathcal{K}}$, which was first solved in \cite{zur}.

Note that $\f \d_d\f\inv=\d_{d'}$ for $d,d'\in\mathcal{T}$, $\f\in\so(\R^3)$ implies
$d=d'$. Thus if $\f:(x,y,z,d,\l)\to(x',y',z',d',\l')$ is a morphism in $\mathcal{K}$ then
$d=d'$ and $\f \d_d\f\inv=\d_d$. Denote by 
$\so_d(\R^3)=\{\f\in\so(\R^3)\mid \f \d_d\f\inv=\d_d \}$ the isotropy subgroup of
$\so(\R^3)$ at $\d_d$. This admits the following description of all morphisms in
$\mathcal{K}$:
$$\Mor((x,y,z,d,\l),(x',y',z',d,\l))=
\{\f\in\so_d(\R^3) \mid (\f(x),\f(y),\f(z))=\!(x',y',z') \} $$
Hence, for every $d\in\mathcal{T}$, we need to find a transversal\footnote{
  A \emph{transversal} for the orbit set of a group action $G\times X\to X$ is a subset of
  $X$ intersecting every $G$-orbit in $X$ in a single point.} 
for the orbit set
$(\R^3\times\R^3\times\R^3)\slash \so_d(\R^3)$ of the action of
$\so_d(\R^3)$ on $\R^3\times\R^3\times\R^3$ by
\begin{equation}\label{gruppverkan}
  \f\cdot(x,y,z)=(\f(x),\f(y),\f(z)).
\end{equation}

Denote by $\iota$ and $\bar{\iota}$ the embeddings $\so(\R^2)\to\so(\R^3)$ given by
$$ \iota(\f)= \left(
\begin{array}{c|c}
  \f&\mbox{$\begin{array}{c}0\\0\end{array}$} \\ \hline
  \rule{0pt}{10pt}\mbox{$\begin{array}{cc}0&0\end{array}$}&1\rule{0pt}{10pt}
\end{array} 
\right) \quad \mbox{and}\quad
\bar{\iota}(\f)= \left(
\begin{array}{c|c}
  1&\mbox{$\begin{array}{cc}0&0\end{array}$}    \\ \hline
  \mbox{$\begin{array}{c}\rule{0pt}{10pt}0\\0\end{array}$} & \f
\end{array} 
\right)$$
respectively.
From \cite{zur}, we take the following description of the isotropy group $\so_d(\R^3)$ for
different $d\in\mathcal{T}$:
\begin{equation*}
  \so_d(\R^3)=
  \begin{cases}
    \so(\R^3) & \mbox{for all} \quad 
    d\in\mathcal{T}_1=\{d\in\mathcal{T} \mid d_1=d_2=d_3\}, \\ 
    \<\iota(\so(\R^2)),\bar{\iota}(-\I_2)\> & \mbox{for all} \quad
    d\in\mathcal{T}_2=\{d\in\mathcal{T} \mid d_1=d_2<d_3\}, \\ 
    \<\bar{\iota}(\so(\R^2)),\iota(-\I_2)\> & \mbox{for all} \quad
    d\in\mathcal{T}_3=\{d\in\mathcal{T} \mid d_1<d_2=d_3\}, \\ 
    \<\iota(-\I_2),\bar{\iota}(-\I_2)\> & \mbox{for all} \quad
    d\in\mathcal{T}_4=\{d\in\mathcal{T} \mid d_1<d_2<d_3\}.
  \end{cases}
\end{equation*}
We denote by $G_i=\so_d(\R^3)$ the isotropy group of any $d\in\mathcal{T}_i$.

Guided by geometric intuition, it is now straightforward to write down normal forms for the
action (\ref{gruppverkan}). For the case $x=0$, we shall reuse the normal forms for the
action
$$G_i\times(\R^3\times\R^3) \to \R^3\times\R^3,\: \f\cdot(x,y)= (\f(x),\f(y)),$$
given in \cite[pp.~17--18]{zur}. If $\mathcal{Q}_i$ is a transversal for the orbit set of
this action, then $\cc_{i0}=\left\{0\right\}\times\mathcal{Q}_i$ is a
transversal for the orbit set of the $G_i$-action on $\{0\}\times\R^3\times\R^3$ induced
by (\ref{gruppverkan}).

Let $N$ be the set of non-negative real numbers, and $P$ the set of positive
ditto. Identifying elements in $\R^3\times\R^3\times\R^3$ with $3\!\times\!3$-matrices, we
have the following normal forms:

\begin{align*}
  \cc_1= \cc_{10}&\cup
  \pmatr{P&\R&\R\\0&0&N\\0&0&0}\cup \pmatr{P&\R&\R\\0&P&\R\\0&0&\R} \\[1ex]
  \cc_2= \cc_{20}&\cup 
  \pmatr{0&0&0\\0&0&N\\P&\R&\R}\cup \pmatr{0&0&\R\\0&P&\R\\P&\R&\R}\cup
  \pmatr{0&0&0\\P&\R&\R\\0&0&N}\cup\pmatr{0&0&P\\P&\R&\R\\0&0&\R}\\  
  &\cup\pmatr{0&0&\R\\P&\R&\R\\0&P&\R}\cup 
  \pmatr{0&P&\R\\P&\R&\R\\0&\R&\R}  
  \cup \pmatr{0&\R&\R\\P&\R&\R\\P&\R&\R} \\[1ex]
  \cc_3=\cc_{30}&\cup
  \pmatr{P&\R&\R\\0&0&N\\0&0&0}\cup \pmatr{P&\R&\R\\0&P&\R\\0&0&\R}\cup
  \pmatr{0&0&\R\\0&0&P\\P&\R&\R}\cup \pmatr{0&0&N\\0&0&0\\P&\R&\R} \\ 
  &\cup\pmatr{0&P&\R\\0&0&\R\\P&\R&\R}\cup  
  \pmatr{0&\R&\R\\0&P&\R\\P&\R&\R}\cup  
  \pmatr{P&\R&\R\\0&\R&\R\\P&\R&\R} \\[1ex]
  \cc_4=\cc_{40}&\cup
  \pmatr{P&\R&\R\\0&0&P\\0&0&\R}\cup \pmatr{P&\R&\R\\0&0&0\\0&0&N} \cup 
  \pmatr{P&\R&\R\\0&P&\R\\0&\R&\R} \cup \pmatr{P&\R&\R\\0&0&\R\\0&P&\R} \\ 
  &\cup\pmatr{0&0&P\\P&\R&\R\\0&0&\R}\cup \pmatr{0&0&0\\P&\R&\R\\0&0&N} \cup
  \pmatr{0&P&\R\\P&\R&\R\\0&\R&\R}\cup \pmatr{0&0&\R\\P&\R&\R\\0&P&\R} \\
  &\cup \pmatr{0&0&P\\0&0&\R\\P&\R&\R}\cup \pmatr{0&0&0\\0&0&N\\P&\R&\R} \cup
  \pmatr{0&P&\R\\0&\R&\R\\P&\R&\R} \cup \pmatr{0&0&\R\\0&P&\R\\P&\R&\R} \\
  &\cup\pmatr{0&\R&\R\\P&\R&\R\\P&\R&\R}\cup \pmatr{P&\R&\R\\0&\R&\R\\P&\R&\R}\cup
  \pmatr{P&\R&\R\\P&\R&\R\\\R&\R&\R}
\end{align*}
Now $\cc_i$ is a transversal for the orbit set of (\ref{gruppverkan}), $i\in\{1,2,3,4\}$,
and
\begin{equation} \label{nform}
\cs=\bigcup_{i=1}^4\,\cc_i\times\mathcal{T}_i\times(\R\minus\{0\})
\end{equation}
a cross-section for the isomorphism classes in $\mathcal{K}$.

Summarising Theorem~\ref{4dim}, Proposition~\ref{dissekv} and Proposition~\ref{konfekv},
we obtain the following classification result:
\begin{cor}\label{classthm}
Every four-dimensional power-commutative division algebra is isomorphic to precisely one
algebra in the set
$$\left\{ \mathcal{F}(\R^3,\xi_x,\eta_{yd},\s_z,\l)
\right\}_{(x,y,z,d,\l)\in\mathcal{S}}.$$
Conversely, all these are power-commutative division algebras of dimension four.
\end{cor}

The multiplication in
the algebra $\mathcal{F}(\R^3,\xi_x,\eta_{yd},\s_z,\l)=\R\times\R^3$ is given by
(\ref{kvadmult}).

\begin{rmk}
  Let $\kappa=(x,y,z,d,\l)\in\mathcal{K}$, and $A=\mathcal{FG}(\kappa)$. An automorphism
  $\kappa\to\kappa$ is an element $\f\in G_d$ fixing the points $x$, $y$ and $z$. Thus it
  is clear that the isomorphism group of $A$ is trivial whenever $x,y,z$ are not on the
  same line through the origin.
  On the other hand, $\Aut(A)\simeq\so(\R^3)$ if and only if $x=y=z=0$ and $d_1=d_2=d_3$,
  that is precisely when $A$ is flexible.
  In this case, $\Aut(A)$ acts transitively on the unit sphere in 
  $\{x\in A\minus\R e \mid x^2\in\R e\}\cup\{0\}\subset A$, which is a subspace supplementary
  to $\R e$.
\end{rmk}

\section{Algebras with a unique non-zero idempotent} \label{sec:enidp}

If a power-commutative division algebra $A$ contains only one non-zero idempotent,
that idempotent is omnipresent in $A$. Hence $A$ is a planar isotope of some quadratic
division algebra $B$. We shall now investigate what conditions a planar isotope 
$A=B_T$ must satisfy to have a unique non-zero idempotent. Assume $\dim A\ge2$. 
Let $\s:\Im B\to \R$ be the linear form and $\l\in\R\minus\{0\}$ the scalar corresponding
to $T$. 

The algebra $A$ has precisely one non-zero idempotent if and only if 
$\idp(\spann\{1_B,v\})\linebreak[2]=\{1_B\}$ for every two-dimensional subalgebra
$\spann\{1_B,v\}\subset A$,  $v\in\Im B\minus\{0\}$.
Let $S=\max\{\s(x) \mid x^2=-1_B\}$. Now multiplication in a arbitrary two-dimensional
subalgebra $U=\spann\{1_B,v\}$ of $A$ is given by
$$(\a1_B + \b v)\circ(\gamma1_B+ \d v)= 
((\a+ s\b)1_B + \l\b v)((\gamma+ s\d)1_B + \l\d v)$$
where $v\in\Im B$ satisfies $v^2=-1_B$, and $s=\s(v)\in[-S,S]$. Choosing the sign of $v$
we may assume $s\in[0,S]$.

Elements in $\idp(U)$ now correspond to non-trivial solutions $(\a,\b)$ to the equation
$(\a1_B + \b v)\circ(\a1_B + \b v)=\a1_B + \b v$. The idempotent $1_B$ corresponds to
$\b=0$, thus we seek solutions with $\b\ne0$.
A calculation gives
$$(\a1_B + \b v)\circ(\a1_B + \b v)=\a1_B + \b v \quad \Leftrightarrow \quad
\begin{cases} 
  \a = (\a+s\b)^2 - (\l\b)^2 & \; (\dagger)\,, \\
  \b = 2\l\b(\a+s\b) & \; (\ddagger)\,,
\end{cases}
$$
and with $\b\ne0$, the identity $(\ddagger)$ is equivalent to $\a=\frac{1}{2\l} - s\b$.
This inserted into $(\dagger)$ gives
\begin{equation} \label{enidp}
\frac{1}{(2\l)^2} -\l^2\b^2 = \frac{1}{2\l} - s\b .
\end{equation}
This is a quadratic equation in $\b$, non-zero solutions of which correspond to
elements in $\idp(U)\minus\{1_B\}$. 

Equation~(\ref{enidp}) can be written as 
$$(\l\b- Ls)^2 = (s^2+1)L^2-L \,,$$ 
where $L=1/(2\l)$.
This equation has solutions in $\b$ if and only if $(s^2+1)L^2-L\ge0$,
that is, if and only if $\l\le(s^2+1)/2$. It has a double root in $\b=0$ if and only if
$s=0$ and $\l=1/2$.

If $\dim A=2$ then $s=S$. If $\dim A >2$ then there are two-dimensional subalgebras such
that $s$ takes any value between $0$ and $S$.
The algebra $A$ has precisely one non-zero idempotent if and only if (\ref{enidp}) has no
non-zero solutions for any possible value of $s$.
In each case, this condition can be formulated as follows.

\begin{prop} \label{enidpsats}
  A planar isotope $A=B_T$ of a quadratic division algebra of dimension at least two has a
  unique non-zero idempotent if and only if either 
  $$\l>(S^2+1)/2\,, \qquad\mbox{or}\qquad  S=0\,,\; \l=1/2\,.$$
\end{prop}

\section{Proof of Theorem~\ref{4dim}} \label{sec:4dim}

In this section we prove that every four-dimensional power-commutative division algebra
contains an omnipresent idempotent.
From here on, let $A$ be a power-commutative division algebras of dimension four.
\begin{lma} \label{centralidp}
  Let $B\subset A$ be a two-dimensional subalgebra, and $e\in A$ an idempotent. If $e$
  commutes with every element in $B$, then $e\in B$. 
\end{lma}

\begin{proof}
We have $[e,x^2]=0$ for all $x\in B$. By the identity~(\ref{kvadkomm}),
$$0=[e,x^2]=[e\bullet x, x]=2[ex,x]$$
so $x$ and $ex$ commute.

Being in itself a division algebra, the subalgebra $B$ contains some non-zero idempotent
$f$ (Lemma~\ref{sdzd}(1)). Let $y\in B$ be any element that is not proportional to $f$.
By Lemma~\ref{kommidpsubalg}, $\spann\{e,f\}$ is a subalgebra of $A$, and since
$[e,y]=[f,y]=0$, we have $\spann\{e,f\}\subset C(y)$. In particular, $[ef,f]=0$ and
$[ef,y]=0$. Now
$$0=[e(f+y),f+y]=[ef,f]+[ef,y]+[ey,f]+[ey,y]= [ey,f].$$

We have proved that the elements $e,f,y,ey$ in $A$ commute with each other.
Note that $e\in B$ if and only if $ey\in B$. Therefore, if $e\not\in B$ then $f,y,e,ey$
are linearly independent, and $A=\spann\{f,y,e,ey\}$. But then $A$ is commutative,
contradicting Hopf's theorem on commutative division algebras. Thus $e\in B$, which was to
be proved.
\end{proof}

\begin{rmk}
  We feel inclined to conjecture that Lemma~\ref{centralidp} holds true also when 
  $\dim A=8$, at least in case the idempotent $e$ is central.
  However, we have not been able to establish a proof that holds in the eight-dimensional
  situation. 
\end{rmk}

An important consequence of Lemma~\ref{centralidp} is that a non-zero central idempotent
$e$ in $A$ must be omnipresent: $A(x)= \spann\{e,x\}$ for every non-idempotent $x\in A$.
Hence, to prove Theorem~\ref{4dim}, it suffices to show that $A$ possesses a central
idempotent.

We shall call an idempotent $e$ in $A$ \emph{exceptional} if $C(e)$ is three-dimensional
and closed under the squaring map $x\mapsto x^2$. The proof of Theorem~\ref{4dim} now proceeds in two steps:
\begin{enumerate} \renewcommand{\labelenumi}{\arabic{enumi}.}
\item The algebra $A$ must contain either a central idempotent or an exceptional
  idempotent.
\item No exceptional idempotent can exist in $A$. 
\end{enumerate}

For the first step, assume that $A$ does not contain a central idempotent. Then, by
Proposition~\ref{eqvillkor}(\ref{Icomm}), there exist $a,b\in\idp(A)$ such that
$[a,b]\ne0$. Lemma~\ref{eminusfkvadrat} implies that $(a-b)^2=a+b - a\bullet b$ is a
multiple of an idempotent $e$, and $e\not\in\spann\{a,b\}$. Since 
$e\in A(a+b)\cap A(a-b)$, we have $\spann\{a,b\}=\spann\{a+b,a-b\}\subset C(e)$. 
By assumption $e$ cannot be central, so $C(e)=\spann\{e,a,b\}=\spann\{a,b, a\bullet b\}$.

\begin{lma} \label{squaring}
The idempotent $e$ is ex\-cep\-tion\-al.
\end{lma}

\begin{proof}
For all $\l,\mu,\nu\in\R$, we have
$$(\l a +\mu b +\nu e)^2 =
\l^2a +\mu^2b +\nu^2e + \l\mu\, a\bullet b + 2\l\nu\, ae +2\mu\nu\, be.$$
By Lemma~\ref{Leinvariant}, $C(e)$ is invariant under multiplication with $e$, so
$ae, be \in C(e)$. Since also $a,b,e,a\bullet b \in C(e)$, this means that 
$(\l a +\mu b +\nu e)^2\in C(e)$. Hence $C(e)$ is closed under squaring, and $e$ is an
exceptional idempotent.
\end{proof}

This gives the first step of the proof of Theorem~\ref{4dim}.
For the second step, suppose that $e\in A$ is an exceptional idempotent (for instance the
one constructed above). 
Note that $\spann\{e,x\}\subset A$ is a subalgebra whenever $x\in C(e)$: by
Lemma~\ref{centralidp}, $\spann\{e,x\}=A(x)$ if $x\in C(e)\minus\idp(A)$, and for
$x\in\idp(A)$, Lemma~\ref{kommidpsubalg} gives the result. 

Let $x\in A\minus C(e)$ be a non-idempotent element. Since $\dim C(e)=3$, the subalgebra
$A(x)$ intersects $C(e)$ in a one-dimensional subspace. Both $C(e)$ and $A(x)$ are closed
under squaring, so $C(e)\cap A(x)$ is closed under squaring. Hence $A(x)$ contains a
unique non-zero idempotent commuting with $e$. 

\begin{lma}
  For every element $x\in A\minus C(e)$ that is not proportional to an idempotent, the
  subspace $A(x)\cap C(e)\subset A$ contains a unique non-zero idempotent, which is
  exceptional.
\end{lma}

\begin{proof}
Since $x$ is not proportional to any idempotent, $\dim A(x)=2$, and $x\not\in C(e)$
implies $\dim A(x)\cap C(e) =1$, thus there is a unique non-zero idempotent $f$ in 
$A(x)\cap C(e)$. By Lemma~\ref{sdzd}(3), there exists an element $y\in A(x)$ such that $y^2=-f$.
Now, for any $\l\in\R$,
$$(\l e+y)^2=\l^2e +\l e\bullet y -f = \l(\l e + y) + \l(e\bullet y - y) - f . $$
By Equation~(\ref{idpkomm}), $[e,e\bullet y - y]=0$, so we have 
$\l(e\bullet y - y) -f \in C(e)$, and $\l(\l e + y)\in\spann\{\l e +y\}$.
This means that $\l(e\bullet y - y) -f= (\l e +y)^2-\l(\l e +y)$ is in 
$C(e)\cap A(\l e + y)$. Since $C(e)\cap A(\l e + y)$ is of dimension one and closed under
squaring, $\l(e\bullet y - y) -f$ must be a scalar multiple (possibly zero) of an
idempotent.

Now we shall show that $e\bullet y -y\in \R f$.
Suppose on the contrary that $e\bullet y -y \,,\, f$ are linearly independent. Then
$\l(e\bullet y - y) -f$ is always non-zero, and proportional to a unique idempotent
$g_\l\in\spann\{e\bullet y - y \,,\, f\}$. Clearly, $g_\l\ne g_\mu$ whenever $\l\ne\mu$,
and thus $\spann\{e\bullet y - y \,,\, f\}$ contains infinitely many idempotents. But this
is impossible, by Lemma~\ref{sdzd}(2).
Hence $e\bullet y -y\in \R f$.

Now $(\l e + y)^2=\l(\l e + y) + \a_\l f \in\spann\{f \,,\, \l e +y\}\subset C(f)$, where
$\a_\l\in\R$. For any $\b\in\R$ and $z\in\spann\{e,y\}$, we have
$(\b f + z)^2 = \b^2f + z^2 +2fz$.
We have shown that $z^2\in C(f)$, and $z\in\spann\{e,y\}\subset C(f)$ implies $fz\in C(f)$.
Consequently, $(\b f + z)^2\in C(f)$. So $C(f)=\spann\{e,f,y\}$ is closed under squaring;
hence $f$ is exceptional.
\end{proof}

Let $x,y\not\in C(e)$ be two element not proportional to any idempotent, and
$A(x)\cap C(e) = \R f$ and $A(y)\cap C(e) = \R g$, for idempotents $f$ and $g$. If $f=g$
then $C(f)\supset\spann\{e,f,x,y\}$, so the latter set must be
linearly dependent (otherwise $f$ would be central). 
This means that $y\in\spann\{e,f,x\}$. 
Conversely, $f\ne g$ whenever $y\not\in\spann\{e,f,x\}$.
If this is the case, and $g\in\spann\{e,f\}$, then $e,f,g$ are three distinct idempotents in
the subalgebra $\spann\{e,f\}\subset A$. Then let $z$ be an element not contained in
$C(e) \cup \spann\{e,f,x\} \cup \spann\{e,f,y\}$ and not proportional to an
idempotent. The idempotent $h\in A(z)\cap C(e)$ given by $z$ is distinct from $e$, $f$ and
$g$. Since there can be no more than three distinct, non-zero idempotents in
$\spann\{e,f\}$, it follows that $e,f,h$ are linearly independent.

The above proves the existence of exceptional idempotents $f,g\in C(e)$ such that
$C(e)=\spann\{e,f,g\}$. Note that $[f,g]\ne0$, otherwise $C(e)$ would be a commutative
subalgebra of $A$. 

Set $C=C(f)\cap C(g)$. Now $C$ has dimension two, and is closed under squaring (since
each of $C(f)$ and $C(g)$ is). So $C$ is a commutative subalgebra, containing
$e$. Thus $C\subset C(e)$. But then 
$$C=C(e)\cap C= C(e)\cap C(f)\cap C(g)=\R e.$$ 
This contradicts $\dim C=2$. We have reached the conclusion that no exceptional idempotent
$e$ can exist. Thus $A$ contains a central idempotent.

\bibliographystyle{plain}
\bibliography{../litt}

\end{document}